\documentclass[10pt,amsfonts, epsfig]{amsart}

\usepackage{amsmath, amscd, amssymb}
\usepackage[frame,cmtip,arrow,matrix,line,graph,curve]{xy}
\usepackage{graphpap, color}
\usepackage[mathscr]{eucal}
\usepackage{mathrsfs}

\numberwithin{equation}{section}

\def\sH{{\mathscr H}}
\def\sO{{\mathscr O}}

\def\sO{\mathscr{O}}

\def\sE{\mathscr{E}}
\def\sF{\mathscr{F}}
\def\sV{\mathscr{V}}

\newcommand{\CC}{\mathbb{C}}

\newcommand{\cal}{\mathcal}

\def\cO{{\cal O}}

\def\mapright#1{\,\smash{\mathop{\lra}\limits^{#1}}\,}

\def\lra{\longrightarrow}

\def\Oplus{\mathop{\oplus}}

\def\begeq{\begin{equation}}
\def\endeq{\end{equation}}
\def\and{\quad{\rm and}\quad}

\def\and{\quad\text{and}\quad}
\def\mapright#1{\,\smash{\mathop{\lra}\limits^{#1}}\,}

\newtheorem{prop}{Proposition}[section]

\newtheorem{lemm}[prop]{Lemma}

\newtheorem{defi}[prop]{Definition}

\def\Ob{\cO b}

\def\blue{\textcolor{blue}}

\def\sO{{\mathscr O}}

\title{Derived Kodaira Spencer map and derived Cosection lemma}
\date{}

\begin{document}

\author{Huai-Liang, Chang}

\maketitle

\begin{abstract}
  The higher tangent vectors and derived Kodaira Spencer maps are defined by K. Behrend and B. Fantechi
  in \cite{Cone}. We collect some properties of them to prove a derived version of the cosection lemma of J. Li
  and Y.H. Kiem (\cite{Jli}) without requiring the perfectness condition. As an application we give a short proof
  of the famous Kodaira's Principle \textit{ambient cohomology annihilates obstruction} (semiregularity).

\end{abstract}

\section{Introduction}


  Let $S$ be a scheme of finite type over $\CC$ and $p\in S$. The cotangent complex $L_S^{\cdot}\in D_{et}(\sO_S)$ (\cite{Illusie, Inf})
  is constructed as a derived functor of Kaehler differential and $T^{\cdot}_S=L_S^{\cdot\vee}$ is the tangent complex of $S$.
  In \cite{Cone} K. Behrend and B. Fantechi  defined the tangent spaces of $S$ at $p$ to be $T^i_{p,S}:=H^{i-1}(T_S^{\cdot}|_p)$, where
  the $T_S^{\cdot}|_p$ is the derived restriction of $T_S^{\cdot}$ at $p$. It is known that $T^1_{p,S}$ is the  Zariski tangent space, $T^2_{p,S}$ is the minimal(intrinsic) obstruction space and  $T^i_{p,S}$($i>1$) called higher tangent spaces. Both $T^1$ and $T^2$ have
  been discussed in many literatures and there is no satisfying interpretation of
 $T^i_{p,S}$ for $i>2$. In this work we study two derived morphisms which map $T^{\cdot}_S$ to perfect complexes: the derived Kodaira Spencer
 map and the derived cosections.\\

 When $Y/S$ is a family of smooth projective varieties there is Kodaira Spencer map from $T^1_{p,S}$ to $H^1(Y_p,T_{Y_p})$.
 By the association $[T\to S] \Rightarrow Y|_T$ the intrinsic moduli problem: deforming $p\in S$, is
 connected to the moduli problem of deforming complex structure of $Y_p$. Hence the intrinsic obstruction space $T^2_{p,S}$
 admits a canonical map to $H^2(Y_p,T_{Y_p})$.  In proposition 6.1 of \cite{Cone}, K. Behrend
 and B. Fantechi extends these two maps to all degrees $T^i_{p,S}\to H^i(Y_p,T_{Y_p})$. They construct a morphism
 in the derived category $\kappa:L_S^{\vee}\to R\pi_*(T_{Y/S})[1]$. The ``derived Kodaira Spencer Map'' $\kappa$ is constructed
 as an application of properties of cotangent complexes. The analogue of $\kappa$ in the bundle deformation problem is the
 Atiyah class constructed by L. Illusie (cf.\cite{Illusie}, \cite{KS}).\\

 The cosection lemma of J. Li and Y-H. Kiem (\cite{Jli}) showed that given a cosection of the obstruction sheaf
 $\xi:\Ob_M\to \sO_M$ over the moduli $M$ with perfect obstruction theory, the reduced intrinsic normal cone $C_p\subset
 \Ob_M|_p$ lies inside the kernel of the map $\xi|_p:\Ob_M|_p\to \sO_M|_p=\CC$. Along with a
 result in \cite{Cone} which said the intrinsic normal cone $C_p$ has its reduced part identical to the collection of curvilinear obstructions, one concludes that the curvilinear obstructions are annihilated by the map $\xi|_p$.\\

 Since many moduli problems are not equipped with perfect obstruction theory, one would expect the cosection lemma
 holds in imperfect case. It is a principle that all moduli problem with moduli space $M$ associates a map from
 the tangent complex of $M$ to the deformation complex $M$. The derived Kodaira Spencer map and Atiyah classes are such
 examples. Hence one expects a version of cosection lemma for the intrinsic deformation problem, which is governed by the
(co)tangent complex $L_M^{\cdot}$. However a cosection like $\xi:\Ob_M^{int}:=H^1(L_M^{\vee})\to \sO_M$ is not
appropriate because the fiber of $\Ob_M^{int}$ at $p\in M$ is not the same as the intrinsic obstruction space
$T^2_{p,S}$ and there is no map from $T^2_{p,S}$ to $\sO_M|_p=\CC$. Hence the vanishing of the cone could not be
discussed. Instead we use a derived cosection $\eta:L^{\vee}_M\to \sO_M[-1]$ and it gives
$H^1(\eta|_p):T^2_{p,S}=H^1(L^{\vee}_M|_p)\to \sO_M|_p=\CC$. So one can discuss the behavior of the intrinsic
cone under this map. We proved the derived cosection lemma with the same conclusion that the reduced intrinsic cone
$C_p$ lies inside the kernel of $H^1(\eta|_p)$.\\

 As an application, we prove the semiregularity under the assumption that local universal family exists.
 The semiregularity is the Kodaira's Principle
 \textit{ambient cohomology annihilates obstruction}: the (curvilinear)obstructions to deformations of $X$ are contained in
 the kernel of
\begin{equation*}
 i:H^2(X,T_X)\rightarrow \Oplus_i Hom (H^{p,q}(X),H^{p-1,q+2}(X)), \quad i_{\xi}(w)=\xi\lrcorner \omega.
\end{equation*}
 It has been proved in a series of paper \cite{Clemens, Manetti, Ran, Manetti1}. Our proof is different from
 theirs in nature.  While the conventional proofs use the differential graded Lie structure of the deformation complex
 $\Oplus H^*(X,T_X)$ at one point, we  combine the derived Kodaira Spencer map with the period map and apply the derived cosection
 lemma. This result also hints that there is an equivalence between the deformation dg-Lie algebra at one point and
 the deformation complex over the moduli.\\

 The author thanks Barbara Fantechi and Si Li for discussion about properties of intrinsic normal cones. There is
 also `higher Kodaira Spencer map''  constructed in \cite{Esnault,canonical} which uses differential operators on
 $S$ instead of the higher tangent vectors or cotangent complexes. The relation between these two is unknown and
 worth studying.\\

\section{Higher tangent vectors}

 We recall the definition of higher tangent vectors of a scheme $S$ at one point $p$ given in \cite{Cone}.
 and give a clarification of their fundamental properties.

\begin{defi}
 Let $S$ be a scheme of finite type over $k$ and $p$ a point in $S$ (closed or not). Let $L_S^{\cdot}\in D_{et}(\sO_S)$
 be the cotangent complex of $S$ defined by \cite{Illusie}. For $i\geqslant 0$ define
$T^i_{p,S}=Ext^{i-1}_S(L_S, k(p))$, as $k(p)$-vector spaces.
\end{defi}

 We have the following proposition. \\
\begin{prop}
$T^i_{p,S}$ satisfy the following (all dual and tensor are derived):\\\\
 (1) $T^i_{p,S}\cong H^{i-1}(L_S^{\vee}\otimes k(p))\cong H^{1-i}(L_S\otimes k(p))^{\vee}$ and $T^i_{p,S}\cong
 H^{-i}(L_{p/S})^{\vee}$ for $i>0$ and $p$ closed point.\\
 (2) $T^1_{p,S}$ is the Zariski tangent space of $S$ at
 $p$. If $T^2_{p,S}=0$ then $S$ is smooth at $p$ and all $T^i_{p,S}=0$ for $i>1$. In general $T^2_{p,S}$ includes
 the obstructions of deforming $p$ in $S$, called intrinsic obstruction space. \\
 (3) Let
 $T=\text{Spec}\,A\to\widehat{T}=\text{Spec}\,A'\to S$ be a small extension at $p\in S$ ($A,A'$ are Artinian
 local $k(p)$ algebras with residue fields $k(p)$ and extension ideal $a$ such that $m_{A'}a=0$).  Then the
 obstruction to extend an element in $Ext^{n-1}_S(L_S, A)$ to be induced from some element in $Ext^{n-1}_S(L_S,
 A')$ lies in $T^{n+1}_{p,S}\otimes_{k(p)} a$. 
 \\
 \end{prop}

\begin{proof}

 (1) Represent $L_S$ by a locally free resolution $G^{\cdot}$ of $\sO_S$ modules then  there is $Ext^{i-1}_S(L_S,
 k(p))=H^{i-1}(Hom(G^{\cdot},k(p))=H^{i-1}(Hom(G^{\cdot},\sO_S)\otimes k(p))=H^{i-1}(L_S^{\vee}\otimes k(p))$.
 Similarly $H^{i-1}(L_S^{\vee}\otimes k(p))\cong H^{i-1}(G^{\cdot\vee}\otimes k(p))=H^{i-1}((G^{\cdot}\otimes
 k(p))^{\vee})=H^{1-i}(L_S\otimes k(p))^{\vee}.$ If $p$ is closed point then the maps $p\to S\to \text{Spec}\,k$
 gives exact triangle $L_S|_p\to L_{p}\to L_{p/S}$ and thus the identity  $T^i_{p,S}\cong
 H^{-i}(L_{p/S})^{\vee}.$\\

(2) The first two statements have been shown by D. Quillen \cite{Quillen}. Let $T=\text{Spec}\,A\to
\bar{T}=\text{Spec}\,A'$ be a square zero extension of ideal sheaf $J$. By (\cite{Illusie}) there is a canonical
element $\omega(g)\in Ext^1(g^*L_S^{\cdot}, J)$ which represents the obstruction of extending a map $g:T\to S$ to
$\bar{T}\to S$.  Assume $T\to \bar{T}$ is a small extension of spectrum of Artinian local rings and $g$ maps the
close point is $T$ to $p\in S$ and denote $\iota:\text{pt}=\text{Spec}\,k(p)\to T$ then  $Ext^1(g^*L_S^{\cdot},\
J)=Ext^1_T(g^*L_S^{\cdot},\iota_*J)=Ext^1_{pt}(\iota^*L_S,J)=Hom_{k(p)}(H^1(L_S\otimes k(p)),J)=H^1(L_S\otimes
k(p))^{\vee}\otimes J=H^1(L_S^{\vee}\otimes k(p))\otimes J=T^2_{p,S}\otimes_{k(p)} J$. In the other words the
deformation of p in $S$ admits an obstruction theory with value in the $k(p)$-vector space $T^2_{p,S}$. \\

(3) The short exact sequence of $\sO_S$ modules: $0\to a\to A'\to A \to 0$ gives exact $$ \cdots\rightarrow
Ext^{n-1}_S(L_S,A') \rightarrow Ext^{n-1}_S(L_S,A) \rightarrow Ext^n_S(L_S,a)\rightarrow \cdots.$$ Since $a$ is
$k(p)$ module by the small extension, there is $Ext^n_S(L_S,a)=Ext^n_{pt}(L_S\otimes k(p),
a)=H^n(L_S^{\vee}\otimes k(p))\otimes_{k(p)} a=T^{n+1}_{p,S}\otimes_{k(p)} a.$

\end{proof}

 Denote $Ext^{n-1}_S(L_S, A)$ by $T^i_{\text{Spec}\,A/ S}$. The above interprets element in $T^{n+1}_{p,S}$ as
 obstructions of extending an element from $T^i_{\text{Spec}\,A/ S}$ to $T^i_{\text{Spec}\,A'/ S}$ for a small
 extension $A'\to A$. This explains in turns that $T^3_{p,S}$ is obstructions of lifting the ''obstruction spaces''
 $T^2_{\text{Spec}\,A/ S}$.\\

 The cotangent complex $L_S^{\cdot}\in D_{et}(\sO_S)$ (\cite{Illusie, Inf}) is constructed locally as
 $\Omega_B\otimes_B A$ where $B\to A$ is a differential graded resolution of $A$ and $\text{Spec}\,A\subset S$
 an affine open chart. The tangent complex $T_S^{\cdot}$ is defined as $L_S^{\cdot\vee}$ and $T_S^{\cdot}|_p$ is equipped with
 a differential graded Lie structure by L. Avramov in \cite{Inf}. There is also a canonical cone inside the second
 tangent space $T^2_{p,S}$. We list these properties here.

\begin{prop}

(1) There exists a canonical differential graded Lie algebra whose complex is $L_S^{\vee}\otimes k(p)[-1]$ and
its cohomlogy $\Oplus_i T^i_{p,S}$ is associated a canonical graded Lie algebra structure.\\
(2) (Quillen's conjecture, Halperin $\&$ Avramov's proof) $T^i(p,S)=0$ for some $i>2 \Leftrightarrow
\text{Spec}\,\sO_{p,S}$ is (local) complete intersection $\Leftrightarrow$   $T^i(p,S)=0$ for all $i>2$. \\
(3) There is a canonical cone inside $T^2_{p,S}$ called intrinsic cone denoted by $C_p$. Its reduced part $C_p^{red}$ is a collection of all curvi-linear obstructions.\\
\end{prop}

\begin{proof}
(1)\cite{Inf}. (2)\cite{conj}. (3) Use Mori's construction (apply the deformation to normal cone to 2.A.11 in
\cite{Moduli} ) to get the cone. For curvi-linear interpretation see proposition 4.7 in \cite {Cone} . (4)
proposition 2 in \cite{Kresch} and section 4 in \cite {Cone}.
\end{proof}

\section{Derived Kodaira Spencer Map}

 From now on we assume the ground field $k$ is equal to $\CC$. The construction
 of the derived Kodaira Spencer map has been given in proposition 6.1 in \cite{Cone}. We recall the construction
 and some properties in this section.\\

 Let $\pi:Y\to S$ be a family of smooth projective varieties of dimension $d$ parametrized by a scheme $S$. Here
$\pi$ is a smooth morphism and $S$ is an arbitrary scheme. We denote $L_S^{\cdot\vee}$ by $T_S$ and
$L_{Y/S}^{\cdot\vee}$ by $T_{Y/S}$. From theory of cotangent complexes (\cite{Illusie}) there is an exact
triangle $L_S|_Y \to L_Y\to L_{Y/S}$ over $Y$  (here $\cdot|_Y$ means $\pi^*(\cdot)$). Taking dual it becomes
$L_{Y/S}^{\vee}\to L_Y^{\vee} \to L_S|_Y^{\vee}$. So there is a map $T_S|_Y=(L_S|_Y)^{\vee}\to
L_{Y/S}^{\vee}[1]=T_{Y/S}[1]$ which gives $\kappa\in$ $Hom_Y(T_S|_Y, T_{Y/S}[1])=Hom_S(T_S,R\pi_*(T_{Y/S}[1]))$.
This is the derived Kodaira Spencer map. It satisfies the following properties.\\

(1) From $T_S\to R\pi_*T_{Y/S}[1]$ there is $T^i_{p,S}=H^{i-1}(T_S|_p)\to H^{i-1}(R\pi_*T_{Y/S}[1]|_p) \to
H^{i}(R\pi_*(T_{Y/S}|p))=H^{i}(Y_p,T_{Y_p})$.  \\

(2) In case $S$ is smooth the total space $Y$ is smooth and the triangle  $L_{Y/S}^{\vee}\to L_Y^{\vee} \to
L_S|_Y^{\vee}$ equals $0\to T_{Y/S}\to T_Y \to T_S\to 0.$ The element
$Hom_Y(T_S|_Y,T_{Y/S}[1])=Ext^1_Y(T_S|_Y,T_{Y/S})$ gives a map $T_S\to R\pi_* T_{Y/S}[1]$ which is equal to
classical Kodaira Spencer for 0-th cohomology, and zero for higher cohomology. And for general $S$, the 0-st
cohomology of the map is the expected Kodaira Spencer map for singular base $S$.\\

(3) The map is functorial in the sense being compatible with base change.\\

  The map $H^2(\kappa_p): T^2_{p,S}\to H^2(Y_p,T_{Y_p})$ has the following interpretation.\\

\begin{lemm}
 The 1-st cohomology of $\kappa$ maps the intrinsic obstruction (obstruction of
lifting $\text{Spec}\,A\to S$ to $\text{Spec}\,A'\to S$ at $p\in S$) to the corresponding obstruction in
$H^2(Y_p,T_{Y_p})$ (obstruction of lifting $Y/\text{Spec}\,A$ to $Y/\text{Spec}\,A'$) described by Kuranishi.
\end{lemm}

\begin{proof}
The case of $S$ equals the universal moduli at $p$ follows from proposition 6.1 in \cite{Cone}. The general case
follows theorem 4.5 in \cite{Cone}.
\end{proof}

 For bundle deformation (over $X\times S$) there is also a derived ``Kodaira Spencer map'' constructed by
 Illusie in \cite{Illusie} and the map is named \textit{Atiyah class}. It is a map from $T_S$ to $R\sH
 om_{\pi_*}(\sE,\sE)$ whose cohomology at $p$ gives $A^i_p:T^i_{p,S}\to Ext^i_X(\sE_p,\sE_p)$ and gives
 obstruction theory if $S$ is the universal moduli. $A^i_p$ are constructed as the composition of the following
 two steps:\\
(1)There is a map from $T^i_{p,S}$ to $Ext^i_{\sO_S}(\CC(p),\CC(p))$ constructed by L. Avramov (theorem 10.2.1
(5) in
\cite{Inf}).\\
(2)By Fourier Mukai transform of the family of sheaves on $X\times S$ there is a functor from $D_{et}(\sO_S)$ to
$D_{et}(\sO_X)$ which induces maps $Ext^i_{\sO_S}(\CC(p),\CC(p))\to Ext^i_X(\sE_p,\sE_p)$.\\

 The composition of this two is the same as the Atiyah class's action at $p$. On the other hand, both the derived Kodaira Spencer
 map and the Atiyah class are maps preserving differential graded Lie algebra structures.


\section{Cosection lemma}

 In \cite{Jli} there is a cosection lemma (2.5) for obstruction theories developed by J. Li and Y.H. Kiem. Starting from a (moduli) space $M$ equipped with perfect obstruction theory and assuming a cosection of the obstruction sheaf  $\sigma : \sO b_{M} \to \sO_{M}$, the  intrinsic cone of $M$ has its reduced part lies inside $\text{ker}\,\sigma$. It is a natural question to ask if one can drop the perfectness for the obstruction theory, for example from the intrinsic obstruction $\sO b_{M}^{\text{int}}=H^1(L_{M}^{\vee})$. \\

 We prove a derived version of the cosection lemma. First we need\\
\begin{lemm}
Assume $M\subset U=\text{Spec} \,\CC[x_1,x_2,...,x_m]$ is defined by the ideal
$I=<f_1,f_2,...,f_n>$. As the zero loci of the section $s=(f_1,f_2,...,f_n)\in\Gamma(U,F=\CC^n\times U)$, $M$ has
a perfect obstruction theory from the kuranishi model in \cite{Virtual}: $\sF^{\cdot}=[\sF^{\vee}|_M
\mapright{\eta} \Omega_U|_{M}]$. By equivalence of the notion of perfect obstruction theory in \cite{Virtual} and
\cite{Cone} there is a complex of locally free sheaf  $L_M^{\cdot}=L^i\to\cdots  \cdots\to L^1 \to L^0 $ and

\begin{eqnarray*}
                   &\sF^{\vee}|_M& \mapright{\eta} \ \ \Omega_U|_{M}\\
                   &\downarrow& \qquad \ \ \ \downarrow\\
\cdots \rightarrow &L^1&           \mapright{\xi}  \ \ L^0.\\
\end{eqnarray*}

In this special case one can pick $L^1=\sF^{\vee}|_M$ and $L^0=\Omega_U|_M$ where the two downwards map are isomorphisms.
\end{lemm}

\begin{proof}
 Take $Kos=$ the Koszul's differential graded resolution of $\sO_M$.
 It maps to the $B^{\cdot}=$ semi-free differential graded resolution of $\sO_M$ which is constructed by
 ``killing homotopies'' in each degree with the $Kos$ to be the first step (see \cite{Inf} chapter 6). Apply
 kaehler differential to the dg-homomorphism $Kos\to B^{\cdot}$ and then restricts (tensor) to $\sO_M$ one gets
 the upper diagram. From this there is $L^0=\Omega_U|_{M}$ and $L^1=\sF^{\vee}|_M$.


\end{proof}

 The derived version of cosection lemma is
\begin{lemm}

Let $M$ be an affine scheme of finite type over $\CC$ and $p$ an arbitrary closed point in $M$.
Assume there is a $G^{\cdot}\in D_{et}(\sO_M)$ such that (1) $G^{\cdot}=\sV^{\vee}[1]$ where $\sV$ is a locally
free sheaf on $M$(2) there is a map $\tau:G^{\cdot}\to L_{M}$. Taking cohomology of $\tau$ there is
$\tau_{-1}=H^{-1}(\tau):\sV^{\vee}=H^{-1}(G^{\cdot})\to H^{-1}(L_M)$ and
$\gamma_{1,p}:H^1(\tau^{\vee}|_p)):T^2_{p,M}=H^1(L_{M}^{\vee}|_p)\to
H^1(G^{\cdot\vee}|_p)=H^{-1}(G^{\cdot})^{\vee}|_p=\sV|_p$. Then the image of curvilinear obstructions under $\gamma_{1,p}$ is zero. \end{lemm}

\begin{proof}

  We reduce the problem to the case with perfect obstruction theory and then apply lemma 2.5 in
  \cite{Jli}. Denote $H^1(L_M^{\vee})$ by $\Ob^{int}_M$ and the corresponding bundle of $\sV$ by $V$. Represent
  $L_M$ by a complex of locally free sheaves by lemma 4.1. One can realize $\tau:G^{\cdot}\to L_M$
  as a map of complexes because of the special form of $G^{\cdot}$. From 4.1 there is a
  surjection $\text{ker}\,\eta=\text{ker}\,\xi\to H^{-1}(L_M)$. Hence the map $\tau_{-1}:\sV^{\vee}\to
  H^{-1}(L_{M})$ factors through some $\tilde{\tau}:\sV^{\vee}\to\text{ker}\,\eta $. This map $\tilde{\tau}$
  could be chosen to be a morphism of complexes
\begin{eqnarray*}
&\sV^{\vee}[1]& \mapright{\tilde{\tau}} \ \sF^{\cdot}\\
&\Vert& \qquad \ \ \downarrow \\
&\sV^{\vee}[1]& \mapright{\tau} \ \, L^{\cdot}
\end{eqnarray*}

Taking dual one gets $\tilde{\tau}^{\vee}: \text{coker}(\eta^{\vee}) \to (\text{ker}\,\eta )^{\vee} \to \sV$. As
$\text{coker}(\eta^{\vee})$ is exactly the obstruction sheaf $\Ob_M$ of the kuranishi model of zero loci of
$s\in\Gamma(U,F=\CC^n\times U)$, there is

\begin{eqnarray*}
 &\Ob_M&  \   \mapright{\tilde{\tau}^{\vee}}\ \ \ \  \sV   \\
&\uparrow \iota& \qquad  \quad         \       \ \, \Vert \\
&\Ob^{int}_M&  \mapright{H^1(\tau^{\vee})} \ \ \sV,
\end{eqnarray*}
and its specialization at any closed point $p\in M$
\begin{eqnarray*}
 &\Ob_M|_p=H^1(\sF^{\cdot\vee}|_p)&  \mapright{\tilde{\tau}^{\vee}|_p}\ \sV|_p\\
&\uparrow &  \qquad     \quad                                           \Vert \\
&T^2_{p,M}=H^1(L^{\cdot\vee}|_p)&     \mapright{\gamma_{1,p}}\          \sV|_p.
\end{eqnarray*}

 There is a normal cone  $C_M$ which comes from the perfect obstruction theory $\sF^{\vee}|_M \mapright{\eta} \ \
 \Omega_U|_{M}$. By lemma 2.5 in \cite{Jli} $C_M|_p$ lies inside the kernel of $\tilde{\tau}^{\vee}|_p$. By
 proposition 2.3 (3) and (4) the cone $C_p$ in  $T^2_{p,M}$ is the same as the cone $C_M|_p$ in $\Ob_M|_p$. Hence
 $\gamma_{1,p}$ maps the cone $C_p^{red}$ to zero by the above diagram. By results in \cite{Cone},  $C_p^{red}$ is
 identical to the collection of curvilinear obstructions of deforming points on $M$ near by $p$.
\end{proof}

\section{Semiregularity}
 Given a smooth compact Kaehler manifold $X$. The differential of the period map at $X$ is
$$i:H^1(X,T_X)\rightarrow \Oplus_i Hom (H^{p,q}(X),H^{p-1,q+1}(X)), \quad i_{\xi}(w)=\xi\lrcorner \omega.$$
 It is known in \cite{Clemens, Manetti, Ran} that obstructions to deformations of $X$ are contained in the kernel of
\begin{equation*}
 i:H^2(X,T_X)\rightarrow \Oplus_i Hom (H^{p,q}(X),H^{p-1,q+2}(X)), \quad i_{\xi}(w)=\xi\lrcorner \omega.
\end{equation*}

The fact is known as Kodaira's Principle \textit{ambient cohomology annihilates obstruction}, under the name of
``semiregularity''. While the proof in \cite{Clemens, Manetti, Ran} doesn't involve properties of period map, M.
Manetti \cite{Manetti1} gives a proof that uses the differential graded Lie algebra structure of $\Oplus_i
H^i(X,T_X)$ and its compatibility with the period map. Assume the existence of an (etale) locally universal
family $Y\to M$ in which $X=Y_0$, we apply the derived Kodaira Spencer map together with the cosection lemma to
give a different and shorter proof. \\

 \begin{proof}
  From $T_{Y/M}\otimes \Omega^p_{Y/M}\to \Omega^{p-1}_{Y/M}\,$ there is
$\,R\pi_*T_{Y/M}\otimes R\pi_*\Omega^p_{Y/M}\to R\pi_*\Omega^{p-1}_{Y/M}$. This gives

$$R\pi_*T_{Y/M}\to R\sH om(R\pi_*\Omega^p_{Y/M},R\pi_*\Omega^{p-1}_{Y/M}).$$

One can shrink $M$ small enough so that the complex

$$K^{\cdot}=R\sH om(R\pi_*\Omega^p_{Y/M},R\pi_*\Omega^{p-1}_{Y/M})$$

 is quasi-isomorphic to its cohomology with zero boundary map

 $$ \cdots\mapright{0}H^i(K^{\cdot})\mapright{0}H^{i+1}(K^{\cdot})\mapright{0}\cdots.$$

Denote $\sV=H^2(K^{\cdot})$ then there is a map

$$R\sH om(R\pi_*\Omega^p_{Y/M},R\pi_*\Omega^{p-1}_{Y/M}) \mapright{j} \sV[-2],$$

such that $H^2(j)$ and $H^2(j|_0)$ are both isomorphisms.
Combine this with the derived Kodaira Spencer map $\kappa: T_M\to R\pi_*T_{Y/M}[1]$, there is
$$T_M\rightarrow R\pi_*T_{Y/M}[1]\rightarrow R\sH om(R\pi_*\Omega^p_{Y/M},R\pi_*\Omega^{p-1}_{Y/M})[1]\rightarrow
\sV[-1].$$

Denote the dual of the complex $\sV[-1]$ by $\,G^{\cdot}\in D_{et}(\sO_M)$. Then the dual of the above
composition is a map $\tau:G^{\cdot}\to L_M$ and $G^{\cdot}$ satisfies the derived cosection lemma. Apply lemma
4.2 for $0\in M$ one has the composition of the map

$$T^2_{0,S}=H^1(L^{\cdot\vee}|_0) \rightarrow  H^2(Y_0,T_{Y_0}) \rightarrow \Oplus_q \sV|_0=Hom
(H^{p,q}(X),H^{p-1,q+2}(X))$$ to be zero. By the lemma 3.1, the image of $C_p$  in $H^2(Y_0,T_{Y_0})$   is the curvilinear obstructions, which vanishes in $Hom (H^{p,q}(X),H^{p-1,q+2}(X))$. This proves the
semiregularity.
 \end{proof}

Remark:
(1)One can prove the semiregularity without using the whole Kodaira Spencer map but instead using a
$[0,1]$ truncation of it. In \cite{KS} the truncated Kodaira Spencer map is constructed and the map $T^2_{p,S}\to
H^2(Y_p,T_{Y_p})$ can be deduced from it. We include the construction of the full Kodaira Spencer map because it
could be useful to study higher tangent vectors. If there is no assumption on the existence of
local universal family one needs to use the stack formulation to deduce the semiregularity.\\
(2) In \cite{Manetti1}, M. Manetti puts a graded Lie algebra and $L_{\infty}$ algebra structure on two sides of
the period map at one point to derive the semiregularity; while here we use the local moduli without algebraic
structures of the deformation complex. The relation is that the information of the moduli space should be encoded
inside the deformation complex at one point with the right algebraic structure. The Maurer Cartan equation is
such an example. There should be more relation between local deformation algebra at $p$ and the moduli space with
the deformation complex (algebra) ($R\pi_*(T_{Y/S})$) to \blue{ALL} orders (Maurer Cartan is of order $\leq 2$).


\begin{thebibliography}{5}
\bibitem{Inf} L. Avramov: \textit{Infinite free resolutions}, in: Six lectures on commutative algebra (Bellaterra, 1996), Progr. Math. 166, Birkhäuser, Basel, 1998; pp. 1-118. \\
\bibitem{conj} L. Avramov: \textit{Locally complete intersection homomorphisms and a conjecture of Quillen on the vanishing of cotangent homology}, Ann. of Math. (2) 150 (1999), 455-487. \\
\bibitem{Cone} K. Behrend, B. Fantechi: The intrinsic normal cone. \textit{Invent. Math.} 128 (1997), no. 1, 45--88.\\
\bibitem{Clemens} H.Clemens: \textit{Geometry of formal Kuranishi theory.} Adv. Math. 198(2005), no. 1, 311-365\\
\bibitem{Esnault} H. Esnault, E. Viehweg: \textit{Higher Kodaira-Spencer classes}, Mathematische Annalen. Volume 299, p 491-527.\\
\bibitem{Manetti1} D Fiorenza, M. Manetti: $L_{\infty}$-algebras, Cartan homotopies and period maps. arxiv: math. AG/0605297.\\
\bibitem{Moduli} D, Huybrechts, M, Lehn: \textit{The Geometry of Moduli Spaces of Sheaves}, Aspects of Mathematics. E, V. 31.\\
\bibitem{KS}  D. Huybrechts, R. P. Thomas: Deformation-obstruction theory for complexes via Atiyah and Kodaira--Spencer classes.  math.AG. arXiv:0805.3527
\bibitem{Illusie} L.Illusie: \textit{Complexe cotangent et deformations I,II.}  Lecture Notes in Mathematics Nos. 239,283. Springer, Berlin, Heidelberg, New York, 1971.\\
\bibitem{Jli} Y.H. Kiem, J. Li:  Gromov-Witten invariants of varieties with holomorphic 2-forms. math.AG/0707.2986\\
\bibitem{Kresch} B. Kim, A. Kresch, T. Pantev:Functoriality in intersection theory and a conjecture of Cox, Katz, and Lee (2001) \\
\bibitem{Virtual}  J. Li,  G. Tian: Virtual moduli cycles and Gromov-Witten invariants of algebraic varieties.  \textit{J. Amer. Math. Soc. } 11  (1998),  no. 1, 119--174.\\
\bibitem{Manetti} M. Manetti: \textit{Cohomological constraint to deformations of compact Kaehler manifolds.} Adv. Math. 186(2004) 125-142; arxiv: math. AG/0105175.\\
\bibitem{Quillen} D. Quillen: \textit{On the (co-)homology of commutative rings}, in: Applications of Categorical Algebra(New York, 1968), Proc. Symp. Pure Math., vol. 17, Amer. Math. Soc., Providence, RI, 1970, pp. 65-87.\\

\bibitem{Ran} Z. Ran: \textit{Universal variations of Hodge structure and Calabi-Yau-Schottky relations.} Invent. Math. 138(1999) 425-449.\\
\bibitem{canonical} Z. Ran: \textit{Canonical Infinitesimal Deformations.}   J. Algebraic Geom.  9  (2000),  no. 1, 43--69.  math.AG/9810041\\

\end{thebibliography}
\end{document}